\documentclass[twoside]{article}

\usepackage{graphicx}
\usepackage{amsmath}
\usepackage{amssymb}
\usepackage{hhline}
\usepackage{epsfig}

\newtheorem{theorem}{Theorem}

\newtheorem{theoremb}{Theorem}
\newtheorem{theoremc}{Theorem}
\newtheorem{theoremd}{Theorem}
\newtheorem{cor}[theoremd]{Corollary}
\newtheorem{dfn}[theoremb]{Definition}

\newtheorem{lem}[theorem]{Lemma}

\newtheorem{rk}[theoremc]{Remark}

\newenvironment{proof}[1][Proof]{\textbf{#1. }}{\qed}
\newcommand\bib[1]{\bibitem[#1]{#1}}

\newcommand\qed{\phantom{\underline{y}}\hfill\hfill$\square$}
\newcommand\comm[1]{}

\newcommand\h{h_{\text{\rm top}}}
\newcommand\hm{H_{\text{\rm mult}}}
\newcommand\hs{H_{\text{\rm sing}}}
\newcommand\e{\epsilon}
\renewcommand\l{\lambda}
\newcommand\op[1]{\mathop{\rm #1}\nolimits}
\newcommand\po{$\!\!\!{\text{\bf.}}$ }
\newcommand\N{{\mathbb N}}
\newcommand\Z{{\mathbb Z}}
\newcommand\R{{\mathbb R}}

\newcommand\Y{{\mathcal Y}}

\newcommand{\weg}[1]{}
\newcommand\ls{\mathop{\overline\lim}\limits_{n \rightarrow \infty}}
\newcommand\ZZ{{\mathcal Z}}

\makeatletter
\renewcommand{\@oddhead}{\hfil Entropy via multiplicity \hfil}
\renewcommand{\@evenhead}{\hfil Boris Kruglikov, Martin Rypdal \hfil}
\makeatother

\begin{document}

%%%%%%%%%%%%%%%%%%%%%%%%%%%%%%%%%%
%%%%%%%%%%%%%%%%%%%%%%%%%%%%%%%%%%

\title{Entropy via multiplicity}
\author{B. Kruglikov  \& M. Rypdal \\ ~ \\
{\small Institute of Mathematics and Statistics}\\
{\small University of Troms\o, N-9037 Troms\o, Norway}\\
{\small Boris.Kruglikov@matnat.uit.no;
Martin.Rypdal@matnat.uit.no} }
\date{}
\maketitle

%%%%%%%%%%%%%%%%%%%%%%%%%%%%%%%%%%
\begin{abstract}
The topological entropy of piecewise affine maps is studied. It is
shown that singularities may contribute to the entropy only if
there is angular expansion and we bound the entropy via the
expansion rates of the map. As a corollary we deduce that
non-expanding conformal piecewise affine maps have zero
topological entropy. We estimate the entropy of piecewise affine
skew-products. Examples of abnormal entropy growth are provided.%
\end{abstract}

%%%%%%%%%%%%%%%%%%%%%%%%%%%%%%%%%%%%%%%%%%%%%%%%%%%%%%%%%%

\section*{Introduction}
For a smooth map $f$ of a compact manifold the Ruelle-Margulis
inequality together with the variational principle \cite{Br,KH}
tells us that the topological entropy of $f$ is bounded by the
maximal sum over positive Lyapunov exponents. For maps with
singularities this result is no longer true and there are examples
of piecewise smooth maps, where the topological entropy exceeds
what can be predicted from the rate of expansion.

In this paper we study the class of piecewise affine maps. It
follows from \cite{B1} that for piecewise affine maps the entropy
is bounded by the rate of expansion and the growth in the
multiplicity of singularities. The latter was shown by J. Buzzi to
be zero for piecewise isometries \cite{B2}, but his proof does not
generalize to non-expanding piecewise affine maps. In fact, we
exhibit an example of a piecewise affine contracting map with
positive topological entropy.

We show that the growth of multiplicity is as effect caused by
angular expansion that can be estimated by the expansion rates of
the map $f$. As a corollary we obtain that for piecewise conformal
maps the topological entropy can be estimated by its expansion
rate as in the smooth compact case. It follows that piecewise
affine non-expanding conformal maps have zero topological entropy.

In the second part of the paper we study the topological entropy
of piecewise affine maps of skew-product type and obtain a formula
which bounds the entropy of the skew products in terms of the
entropy and multiplicity growth of the factors. The estimate
includes a term which indicates that the entropy of a skew-product
system may be greater than the sum of the maximal entropy of its
factors and we give an example where this is realized.

Our main results have several corollaries for which one can
calculate the topological entropy of various classes of piecewise
affine maps.

%%%%%%%%%%%%%%%%%%%%%%%%%%%%%%%
%% DEFINITIONS AND MAIN RESULTS
%%%%%%%%%%%%%%%%%%%%%%%%%%%%%%%

\section{Definitions and main results}

%%%%%%%%%%%%%%%%%%%%%%%%%%%%%%%
%% Piecewise affine maps
%%%%%%%%%%%%%%%%%%%%%%%%%%%%%%%

\subsection{Piecewise affine maps and topological entropy}
\begin{dfn}\po
We say that $(X,\ZZ,f)$ is a piecewise affine map if
\begin{enumerate}
\item $X \subset \R^n$
\item $\ZZ=\{Z\}$ is a finite collection of open, pairwise
disjoined polytopes such that $X':=\cup_{Z \in \ZZ} Z$ is dense in
$X$.
\item $f_Z:=f|_Z:Z \rightarrow X$ is affine for each $Z \in \ZZ$
\end{enumerate}
\end{dfn}

The maps $f_Z$ are called the affine components of the map $f$.
The linear part of $f_Z$ is denoted $f'_Z$, and if $x \in Z$ we
denote $d_x f=f_Z'$. Let $\op{PAff}(X;X)$ be the set of piecewise
affine maps on $X$ and let $U_f=X'\cap f^{-1}(X') \cap f^{-2}(X')
\cap \dots$ be the set of points in $X$ with well-defined infinite
orbits. Let $\ZZ^n$ be the continuity partition of the piecewise
affine map $f^n$. We always assume $X$ to be compact.

Since the maps we consider in this paper have singularities
(consult \cite{KS,ST}) we must define what we mean by topological
entropy.

Let $U_n=\cap_{k=0}^{n-1} f^{-k}(X')$, then $U_f=\cap_{n \geq 0}
U_n$. Take a metric $d$ defining the standard topology on $X$. Let
$d_n^f=\max_{0 \leq k <n} (f^k)^* d$, and define
$S(d_n^f,\epsilon)$ to be the minimal number of
$(d_n^f,\epsilon)$-balls needed to cover $U_n$. Define
$$\h(f)=\lim_{\epsilon \rightarrow 0}\ls \frac1n \log
S(d_n^f,\epsilon)\,.$$ This number is independent of the choice of
metric on $X$ and is finite because it is bounded by $d \cdot
\sup_x \log(\|d_x f\|+|\ZZ|)$. It equals the $(n,\e)$-entropy
$\h(f|_{U_f})$, which coincides with the (upper=lower) capacity
entropy $Ch_{U_f}(f)$ \cite{Pe}. This bounds the topological
entropy $h_{U_f}(f)$ of non-compact subsets by Pesin and Pitskel'
\cite{Pe}, so that we have $h_{U_f}(f)\leq \h(f)$. Since we
estimate $\h(f)$ from above, this bounds the other entropy too.

\begin{rk}\po\label{rk1}
We observe that even in the presence of singularities the property
$\h(f^T)=T\h(f^T)$ holds for $T\in \N$. The proof uses the fact
that for all $\e >0$ there is a number $\delta(\epsilon)>0$ such
that $B_d(x,\delta(\e)) \subset B_{d^f_T}(x,\e)\,\forall x \in
X'$, cf.\,\cite{KH}.
\end{rk}

One can also measure the orbit growth of a piecewise affine map
through the growth of continuity domains. The singularity entropy
of $f$ is
$$\hs(f)=\ls \frac1n \log  |\ZZ^n|\,.$$
For non-expanding piecewise affine maps $\h(f) \leq \hs(f)$. If
$f_Z(x)\ne f_{Z'}(x)$ for all $x \in \partial Z \cap
\partial Z'$ with $Z \neq Z'$, then
$\hs(f) \leq \h(f)$ \cite{R}.

%%%%%%%%%%%%%%%%%%%%%%%%%%%%%%%
%% Lyapunov numbers and multiplicity entropy
%%%%%%%%%%%%%%%%%%%%%%%%%%%%%%%

\subsection{Expansion rates and multiplicity entropy}
\begin{dfn}\po
The multiplicity of the partition $\ZZ^n$ at a point $a \in X$ is
$\mbox{\em mult}(\ZZ^n,a)= |\{Z \in \ZZ^n\,|\, \bar{Z} \ni a\}|$,
and the multiplicity of $\ZZ^n$ is $\mbox{\em
mult}(\ZZ^n)=\sup_{a\in X} \mbox{\em mult}(\ZZ^n,a)$. The
multiplicity entropy \cite{B1} of $f$ is defined as
$$
\hm(f)=\ls \frac1n \log \mbox{\em mult}(\ZZ^n)\,.
$$
\end{dfn}

\begin{dfn}\po
For $f \in \op{PAff}(X,X)$ we define
$$
\lambda^+(f)= \ls \sup_{x \in U_n} \frac1n \max_{0 \leq k \leq d}
\log \|\Lambda^k d_x f^n\|\,.
$$
We also let (the second quantity can equal $-\infty$ for
non-invertible maps)
 $$
\lambda_{\op{max}}(f)=\ls \sup_{x \in U_n} \frac1n \log \|d_x
f^n\| \text{ and } \lambda_{\op{min}}(f)=
 %\mathop{\underline\lim}\limits_{n\to \infty}
 %\inf_{(x,v) \in STU_n} \frac1n \log \|d_x f^n(v)\|=
-\ls\sup_{x \in U_n} \frac1n \log \|(d_x f^n)^{-1}\|\,.
 $$
\end{dfn}

\begin{theorem}\label{th1}\po
For any $f \in \op{PAff}(X;X)$ it holds:
$$
\h(f) \leq \lambda^+(f)+\hm(f)\,.
$$
\end{theorem}

This result is basically due to Buzzi \cite{B1}. However, he only
proves it for a special class of strictly expanding maps, and he
considers the entropy of coding $\hs(f)$ instead of $\h(f)$. Hence
we modify his proof.

%%%%%%%%%%%%%%%%%%%%%%%%%%%%%%%
%% Angular expansion and sphericalizations
%%%%%%%%%%%%%%%%%%%%%%%%%%%%%%%

\subsection{The spherization and angular expansions}\label{S13}
 \comm{
For smooth maps of compact spaces the topological entropy is
bounded by the number $\lambda^+(f)$. This is not true for
piecewise affine maps and a counterexample is given in Example 2.
In this example the mechanism for growth of entropy is the angular
expansion which generates a point with exponential growth of
multiplicity. To give estimates for the topological entropy of
piecewise affine maps we must therefore introduce the maximal
rates of angular expansion, i.e the expansion rates of the
spherization of $f$.
 }
For any $k$-dimensional submanifold $N^k \subset X$ and $x \in N$
we define the spherical bundle $STN=\{v \in TN:\|v\|=1\}$, where
$\|\cdot\|$ is the Euclidian norm on every $T_xN \subset \R^d$.
Let $f$ be non-degenerate, i.e. each affine component is
non-degenerate. The spherization of $f$ is defined to be the
piecewise smooth map $d_x^{(s)}\!f:S T_x X \to S T_{f(x)}X$ given
at $x\in X'$ by the formula
$$
d_x^{(s)}\!f(v)=\frac{d_x f(v)}{\|d_x f(v)\|}\,.
$$
For $x\in \op{Sing}(X)\stackrel{\text{def}}=X\setminus X'$ and
$v\not\in T_x\op{Sing}(X)$ (the tangent cone) we define
$d_x^{(s)}\!f(v)=\lim\limits_{\e\to+0}d_{x+\e v}^{(s)} f(v)$. For
other $(x,v)\in STX$ the map is not defined.

The angular expansion of $f$ is exactly the expansion in its
spherization.

\begin{dfn}\po
For a non-degenerate map $f \in \op{PAff}(X,X)$ and $i<d$ we
define
$$
\rho_i(f)=\ls \frac1n \sup_{x \in U_n} \max_{0 \leq k \leq i}
\sup_{v \in S^{d-1}} \log \| \Lambda^k d_vd_x^{(s)}\!f^n \|.
$$
\end{dfn}

\noindent Note that
 $$
\rho_i(f)= \ls \frac1n \,\sup_{N^i \subset X}\, \sup_{x\in U_n\cap
N^i}\, \sup_{v \in ST_x N^i}\, \log\|\Lambda\, d_v
d^{(s)}_x\!f^n\|\,.
 $$

The number $\rho_{1}(f)$ measures the maximal exponential rate
with which angles can increase under the map $f$. The numbers
$\rho_i(f)$ for $i<d$ measures the maximal rate of expansion of
the restrictions to $(i-1)$-dimensional spheres. Clearly
$\rho_0(f)=0$ for any $f \in \op{PAff}(X,X)$, and if $f$ is
conformal, i.e. all the affine components of $f$ are conformal,
then $\rho_i(f)=0$ for all $i$.

\begin{theorem} \label{th2} \po
$\hm(f) \leq \sum_{i=1}^{d-1} \rho_i(f)$ for any non-degenerate $f
\in \op{PAff}(X;X)$.
\end{theorem}

\noindent The following corollaries are direct consequences of
Theorem \ref{th2}.

\begin{cor}\po\label{cor1}
If $f \in \op{PAff}(X;X)$ is conformal, then $\h(f)\leq
\lambda^+(f)$.
\end{cor}

We say that a piecewise affine map is non-expanding if all its
affine components are non-expanding, i.e. the eigenvalues of the
linear part of each affine component have absolute values not
exceeding $1$.

\begin{cor}\po\label{cor2}
If $f \in \op{PAff}(X;X)$ is conformal non-expanding,
 %(for instance a piecewise isometry),
then $\h(f)=0$.
\end{cor}

It is shown in \S\ref{S31.} that
$$
\rho_i(f) \leq \ls \sup_{x \in U_f} \max_{0 \leq k \leq i} \frac1n
\log \|\Lambda^k d_x f^n \|-i \lambda_{\op{min}}(f) \leq i
\Big{(}\lambda_{\op{max}}(f)-\lambda_{\op{min}}(f) \Big{)} \,,
$$
This gives the following estimate:
\begin{cor}\po\label{cor3}
For a non-degenerate map $f \in \op{PAff}(X;X)$ it holds:
$$
\hm(f)\leq \frac{d(d-1)}{2} \big{(}\lambda_{\op{max}}(f)
-\lambda_{\op{min}}(f)\big{)}\,.
$$
\end{cor}

\noindent Hence we see that the topological entropy of a
non-degenerate piecewise affine map $f$ can be bounded using only
its expansion rates. In fact, we have
$$
\h(f) \leq   \lambda^+(f) +\frac{d(d-1)}{2}
\big{(}\lambda_{\op{max}}(f)-\lambda_{\op{min}}(f)\big{)}\,.
$$

Let us call $f \in \op{PAff}(X;X)$ asymptotically conformal if
$\lambda_{\op{max}}(f)=\lambda_{\op{min}}(f)$.

\begin{cor}\po\label{cor4}
For asymptotically conformal $f \in \op{PAff}(X;X)$: $\h(f) \leq
\lambda^+(f)$.
\end{cor}

 \begin{rk}\po
It is essential that in our definition of piece-wise affine maps
we consider a finite number of continuity domains. With countable
number of domains (this is related to countable Markov chains) the
above theorems become wrong. In fact, according to \cite{AOW}
every aperiodic measure preserving transformation can be
represented as an interval exchange with countable number of
intervals. In particular, there exist piece-wise isometries with
infinite number of continuity domains, which have positive
topological and metric entropies.
 \end{rk}

%%%%%%%%%%%%%%%%%%%%%%%%%%%%%%%
%% Piecewise affine skew products
%%%%%%%%%%%%%%%%%%%%%%%%%%%%%%%

\subsection{Piecewise affine skew products}
\begin{dfn}\po
We say that $S\tilde\times T \in \op{PAff}(X \times Y,X \times Y)$
is a piecewise affine skew product if it has the form
$f(x,y)=(S(x),T_x(y))$ for some $T_x \in \op{PAff}(Y,Y)$.
\end{dfn}

\noindent The following results hold for piecewise affine skew
products:

\begin{theorem} \po \label{th3}
If $S\tilde\times T$ is a piecewise affine skew product, then
$$
\h(S) \leq \h(S\tilde\times T) \leq \h(S) + \hm(S)+\sup_{{\bf x}}
\big{(} \lambda^+(T_{\bf x})+ \hm(T_{\bf x}) \big{)}\,.
$$
where ${\bf x}=(x_0,x_1,\dots)$ is an orbit of $S$ in $X$ and
$T_{\bf x}$ is the dynamics along this orbit, i.e. $T_{\bf
x}^n=T_{x_{n-1}}\circ \dots \circ T_{x_0}$ (see \S\ref{S23} for
details).
\end{theorem}

From Theorem \ref{th3} we can deduce several simple corollaries:

\begin{cor}\po\label{cor5}
If $\dim(X)=1$ and $T_x\in \op{PAff}(Y;Y)$ are non-expanding and
conformal for all $x \in X'$, then $\h(S \tilde\times T) = \h(S)$.
\end{cor}

\begin{cor}\po \label{cor6}
Let $X=[0,1]^d$ and $A \in \op{PAff}(X;X)$ be defined by $x
\mapsto Ax \mod \Z^d$ for some $A \in\op{GL}_d(\R)$. If $T_x \in
\op{PAff}(Y;Y)$ are non-expanding and conformal for all $x \in
X'$, then
$$\h(A \tilde\times T) = \h(A) = \log \, \text{{\em Jac}}^+ A \,,$$
where
$$
\op{Jac}^+A=\prod_{\lambda \in \op{Sp}(A)} \max \{|\lambda|,1
\}\,.
$$
\end{cor}

\begin{cor}\po \label{cor7}
Let $\Sigma_N^+=\{1,\dots,N\}^{\Z_{\geq 0}}$ and $\sigma_N^+$ be
the right shift on $\Sigma_N^+$. Take $T_1,\dots,T_N \in
\op{PAff}(Y;Y)$ to be non-expanding and conformal. Then for the
map $f:\Sigma_N^+ \times Y \rightarrow \Sigma_N^+ \times Y$,
$({\bf t},y) \mapsto (\sigma_N^+ {\bf t}, T_{t_0}(y))$, we have:
$\h(f)=\log N$.
\end{cor}

\begin{rk}\po
The class of piecewise affine skew-products of the form
$$\sigma_N^+ \tilde\times
T:\Sigma_N^+ \times Y \rightarrow \Sigma_N^+ \times Y
$$
is physically relevant and appears in the Zhang sandpile model of
Self-Organized Criticality \cite{BCK}. The maps $T_i$ correspond
to the avalanches and are contracting \cite{KR1}, though not
conformal. So we can get an estimate for the entropy.
\end{rk}

In general we cannot ensure existence of invariant measures for
piecewise affine maps. In fact, there are examples with no
invariant measure. However we can give an estimate for the metric
entropy whenever such a measure exists.
 % The following result is based on the Abramov-Rokhlin formula.

 \begin{theorem} \po
Let $S\tilde\times T$ be a piecewise affine skew product and $T_x
\in \op{PAff}(Y;Y) $ be non-expanding for all $x \in X'$. If $\mu$
is a $S\tilde\times T$-invariant Borel probability measure on $X
\times Y$, then
$$
h_{\pi_*\mu}(S) \leq h_\mu (S\tilde\times T) \leq h_{\pi_*\mu}(S)
+ \hm(S)\,,
$$
where $\pi:X \times Y \rightarrow X$ is the projection to $X$.
 \end{theorem}

\begin{cor}\po\label{cor8}
Let $A \tilde\times T$ be as in Corollary \ref{cor6} with $A$
expanding. If $\mu$ is a measure of maximal entropy for $S
\tilde\times T$ on $X \times Y$, then $\pi_*\mu$ is absolutely
continuous with respect to the Lebesgue measure on $X$.
\end{cor}

%%%%%%%%%%%%%%%%%%%%%%%%%%%%%%%
%% PROOFS OF THE THEOREMS
%%%%%%%%%%%%%%%%%%%%%%%%%%%%%%%

\section{Proof of the theorems}

%%%%%%%%%%%%%%%%%%%%%%%%%%%%%%%
%% PROOF OF THEOREM 1
%%%%%%%%%%%%%%%%%%%%%%%%%%%%%%%

\subsection{Proof of Theorem 1}
Fix $\epsilon >0$ and let $T=T(\epsilon) \in \N$ be such that $
\text{mult}({\mathcal Z}^n)\leq
\exp((H_{\text{mult}}(f)+\epsilon)n)$ for all $n \geq T$ and
$$
\forall x \in U,n \geq T,0 \leq k \leq d:\,\,\|\Lambda^k d_x f^n\|
\leq \exp((\lambda^+(f)+\epsilon)n)
$$
We suppose that $\sqrt{d}+1 \leq \exp(\epsilon T/d)$. Take
$r=r(\epsilon)$ to be compatible with the partition ${\mathcal
Z}^T$, i.e. any $r$-ball intersects maximally
$\op{mult}(\mathcal{Z}^T)$ partition elements.

We will prove inductively on $l$ that each $Z \in \ZZ^{lT}$ can be
covered by a family $Q_Z=\{W\}$ satisfying the following
properties:
\begin{enumerate}
\item $\sum_{Z \in \ZZ^{lT}} \text{card }Q_Z
\leq C_0 \exp((\lambda^+(f)+H_{\text{mult}}(f)+3\e)lT)$
\item $\text{diam}(f^{lT}(W)) \leq r$.
\item $\forall x,y\in W:\,\,d_l^{f^T}\!\!(x,y)<\e$ and
$d_{lT}^f(x,y)<\delta(\e)$ with $\delta(\e)$ from Remark
\ref{rk1}.
\end{enumerate}

\noindent The base of induction $l=0$ is obvious and $C_0 \leq
|\ZZ| (\op{diam} X/\min \{r,\e\})^d$.

Take a partition element $W \in Q_Z$ that is used to cover the set
$Z \in \ZZ^{lT}$. By the induction hypothesis it can be continued
to cover an element of length $\ZZ^{(l+1)T}$ in at most
$\text{mult}({\mathcal Z}^T)$ ways. So to cover the cylinders $Z
\in \ZZ^{(l+1)T}$ we make a division of $W$:
$$
W=\bigcup_{i=1}^\gamma W'_i\,\,,\gamma \leq \text{mult}({\mathcal
Z}^T)\,.
$$
Let $ W''_i=f^{Tl-1}(W'_i)$ and $W'''_i=f^T(W''_i)$. By the
assumption $\text{diam}(W'')\le r$, but the set $W'''$ may have
greater diameter than $r$. Thus  we need to divide the sets
$W'''_i$ and pull this refinement back to the partition of sets
$W'_i$.

The image $f^T(W''_i)$ is the image of one affine component of
$f^T$. Let $L^T$ denote the linear part of this affine component.
We can assume that $L^T$ is symmetric and take $\{e_k\}$ to be a
basis of eigenvectors corresponding to eigenvalues
$\l_1^T,\dots,\l_d^T$. Let $\{v_k\}$ be a basis in the vector
subspace corresponding to $W'''_i$. We can choose this basis to be
orthonormal and triangular with respect to $\{e_k\}$. Divide
$W'''_i$ by the hyperplanes
 $$
\psi_j(x)\stackrel{\text{def}}=\langle v_j,x\rangle=p \frac{\min
\{ r,\epsilon\}}{\sqrt{d}},\quad p\in \Z,\quad j=1,\dots,d\,.
 $$
This defines cells $\tilde W$ of diameter less than
$\min\{r,\epsilon\}$. Since $\psi_j(W''')=\psi_j(L^T(W''))$ has
$\op{diam}\le|\l_i^T| \min\{r,\epsilon\}$, the number of cells
$\tilde W$ needed to cover $W'''$ is less than or equal to
 \begin{eqnarray*}
(\sqrt{d}+1)^d|\l_1^T|^+ \dots |\l_d^T|^+ &\leq& (\sqrt{d}+1)^d
\sup_{x \in U_T}\max_{1\leq k \leq d}\|\Lambda^k d_x f^T\|  \\
&\leq& \exp \big{(}(\lambda^+(f) + 2 \e)T \big{)}\,,
  \end{eqnarray*}
where $|\lambda|^+=\max\{|\l|,1\}$.

Therefore the total cardinality of the new partition is less than
or equal to
 \begin{multline*}
\op{mult}({\mathcal Z}^T)\exp \big{(}(\lambda^+(f)+2\e)T\big{)}
\exp \big{(}(\lambda^+(f)+H_{\text{mult}}(f)+3\e)lT \big{)}  \\
\leq \exp \big{(}(\lambda^+(f)+H_{\op{mult}}(f)+3\e)(l+1)T
\big{)}\,.
 \end{multline*}
The elements of the partition $Q_Z$ have diameter less than $\e$
in the metric $d_l^{f^T}$. By Remark \ref{rk1} each partition
element has diameter less than some number $\delta(\e)>0$ in the
metric $d_{lT}^f$. This proves the statement.

%%%%%%%%%%%%%%%%%%%%%%%%%%%%%%%
%% PROOF OF THEOREM 2
%%%%%%%%%%%%%%%%%%%%%%%%%%%%%%%

\subsection{Proof of Theorem 2}
Define the bundles
 $$
S^{(k)}TX=\{(x,v_1,\dots,v_k)\,|\,x \in X, v_1 \in ST_x X,  v_i
\in ST_x X   \cap \langle v_1,\dots, v_{i-1} \rangle^\bot\}.
 $$
They form the spherical towers:
 $$
S^{(d)}TX\stackrel{\pi_{d}}\longrightarrow
S^{(d-1)}TX\stackrel{\pi_{d-1}}\longrightarrow
\dots\stackrel{\pi_2}\longrightarrow
S^{(1)}TX=STX\stackrel{\pi_1}\longrightarrow X
 $$
with fibers $S^0,S^1,\dots,S^{d-1}$ respectively.

The spherization $d^{(s)}\!f:STX\to STX$ induces the maps
$S^{(k)}\!f:S^{(k)}TX \rightarrow S^{(k)}TX$. Although defined on
$S^{(k)}TX'$ they extend over the strata of $\op{Sing}(X)$ as in
\S\ref{S13} (modulo spherization this corresponds to the
differential $f_{x,v_1,\dots,v_{k-1}}(v_k)$ of Tsujii and Buzzi
\cite{T,B2}):
 $$
S^{(k)}\!f(x,v_1,\dots,v_k)=\lim_{\e_1\to+0}\dots\lim_{\e_k\to+0}
d_{x+\e_1v_1+\dots+\e_kv_k}^{(s)}f(x,v_1,\dots,v_k)
 $$
(the r.h.s. operator is applied to each vector $v_i$
successively). In particular, the map $S^{(d)}\!f$, though has
singularities, is defined everywhere on $S^{(d)}TX$.

Let $\ZZ^n_{(x,v_1,\dots,v_{k-1})}(S^{(k)}\!f|S^{(k-1)}\!f)$ be
the continuity partition of the piecewise smooth map
$(S^{(k)}\!f)^n$ restricted to the fiber
$\pi_{k}^{-1}(x,v_1,\dots,v_{k-1})$. Define
 $$
\hs(S^{(k)}\!f|S^{(k-1)}\!f)= \ls \sup_{(x,v_1,\dots,v_{k-1})}
\frac1n \log
|\ZZ^n_{(x,v_1,\dots,v_{k-1})}(S^{(k)}\!f|S^{(k-1)}\!f)|\,.
 $$

Clearly
$\op{mult}(\ZZ^n_{(x,v_1,\dots,v_{k-1})}(S^{(k)}\!f|S^{(k-1)}\!f),v_k)
=|\ZZ^n_{(x,v_1,\dots,v_k)}(S^{(k+1)}\!f|S^{(k)}\!f)|$. So
$\hm(S^{(k)}\!f|S^{(k-1)}\!f)= \hs(S^{(k+1)}\!f|S^{(k)}\!f)$.

In particular, we have: $\hm(f)= \hs(S^{(1)}\!f|f)$.
 \comm{
Let $\ZZ(S^{(1)}\!f)$ be the continuity partition of the piecewise
smooth map $S^{(1)}\!f$, and let $\ZZ^n(S^{(1)}f)$ be the
continuity partition of its $n$-th power. For $x \in X$ we define
$$
\ZZ^{(s,1),n}_x=\{ Z \in \ZZ^n(S^{(1)}f)\,|\,\bar{Z} \cap S T_x X
\neq \emptyset\}\,,
$$
and let
$$
\hs(S^{(1)} f|f)=  \ls \sup_{x \in X} \frac1n \log
|\ZZ_{x}^{(s,1),n}|\,.
$$
Clearly $\op{mult}(\ZZ^n,x)=|\ZZ_{x}^{(s,1),n}|$, so $\hm(f)=
\hs(S^{(1)}f|f)$.
 }
Applying the arguments from the proof of Theorem 1 we obtain:
 $$
\hs(S^{(1)}f|f) \leq \rho_{d-1}(f)+ \hm(S^{(1)}\!f|f)\,.
 $$
Doing its once more for $\hm(S^{(1)}\!f|f)=
\hs(S^{(2)}\!f|S^{(1)}\!f)$ we get:
 $$
\hs(S^{(2)}f|f) \leq \rho_{d-2}(f)+ \hm(S^{(2)}\!f|f)\,.
 $$
Applying the same argument $d-1$ times yields: $\hm(f) \leq
\rho_{d-1}(f) +\dots + \rho_1(f)$.

%%%%%%%%%%%%%%%%%%%%%%%%%%%%%%%
%% PROOF OF THEOREM 3
%%%%%%%%%%%%%%%%%%%%%%%%%%%%%%%

\subsection{Proof of Theorem 3}\label{S23}
The inequality $\h(S) \leq \h(S \tilde\times T)$ is obvious since
$S$ is a quotient of $S \tilde\times T$. We will now prove the
upper bound for $\h(S \tilde\times T)$. Denote $d_Y=\dim Y$.

Let $\Gamma$ be the set of all sequences
$((x_0,z_0),(x_1,z_1),(x_2,z_2),\dots) \in (X\times
\R^{d_Y})^{\Z_{\geq 0}}$ satisfying $S(x_i+\epsilon z_i)
\rightarrow x_{i+1}$ as $\epsilon \rightarrow+0$. Define
$T_{(x,z)}=\lim_{\e\to+0} T_{x+\e z}$ and let $\Y^{n,{\bf x}}$ be
the collection of non-empty sets
 $$
Y_{(x_0,z_0)} \cap T_{(x_0,z_0)}^{-1}(Y_{x_1,z_1}) \cap \dots \cap
(T_{(x_{n-1},z_{n-1})}\circ \dots
T_{(x_0,z_0)})^{-1}(Y_{(x_n,z_n)})\,,
 $$
for ${\bf x}=((x_0,z_0),(x_1,z_1),(x_2,z_2),\dots) \in \Gamma$,
where the sets $Y_{(x_i,z_i)}$ are the continuity domains of the
maps $T_{(x_i,z_i)}$. Let
 $$
P^n(x)=\{Y \in \Y^{n,{\bf x}}\,|\,{\bf
x}=((x_0,z_0),(x_1,z_1),\dots)\in\Gamma,x_0=x\}
 $$
be the continuity partition iterated along all possible $S$-orbits
starting from $x \in X$. Clearly
$$
|P^n(x)|\leq \text{mult}(S^n,x) \sup_{{\bf x}\in \Gamma} |
\Y^{n,{\bf x}}|\,.
$$
To simplify notations we will write elements of $\Gamma$ as ${\bf
x}=(x_0,x_1,x_2,\dots)$, where $x_i$ consists of a point in $X$
and a vector in $\R^{d_Y}$. If $x_i \in X'$ the vector $z_i$ is
not essential. Denote $T_{\bf x}^n=T_{x_{n-1}} \circ \dots \circ
T_{x_0}$ for ${\bf x}=(x_0,x_1,x_2,\dots) \in \Gamma$ and let
 $$
\hm(T_{\bf x})=\ls \frac1n \log \text{mult}(T_{\bf x}^n)\,.
 $$

\begin{lem}\po \label{lem4}
For a piecewise affine skew product $S \tilde\times T$ with $T_x$
non-expanding for all $x \in X'$, it holds:
$$\h(S \tilde\times T) \leq \h(S) + \hs^{\text{\em fiber}}(T|
S)\,,$$ where
$$\hs^{\text{\em fiber}}(T
|S) = \sup_{x \in X} \ls \frac1n \log | P^n(x)|\,.
$$
\end{lem}

\begin{rk}\po\label{rk5}
The statement of the lemma is similar to Bowen's Theorem 17
\cite{Bo}, but the direct generalization fails, see Example 2.
\end{rk}

\begin{proof}
Let $\epsilon>0$ be arbitrary. Denote $a=\hs^{\text{fiber}}(T|
S)$, and fix $\alpha>0$ and $m_\alpha=[1/\alpha] \in \N$. For all
$x \in X$ we let
 $$
n_\alpha(x)=\min \{n \geq m_\alpha\,|\, \frac1n \log | P^n(x)|
\leq a+\alpha\}\,.
 $$
Since the function $\op{mult}(\mathcal{Z},x)$ is upper
semi-continuous, the same is true for the function $n_\alpha(x)$.
So $n_\alpha:=\sup_{x \in X} n_\alpha(x)$ is finite. Then we have
$0 < m_\alpha \leq n_\alpha(x) \leq n_\alpha < +\infty$ for all $x
\in X$.

Let $r_\alpha>0$ be compatible with all the partitions $\ZZ^{n}$,
$m_\alpha\le n\le n_\alpha$, i.e. any ball of radius $r_\alpha$
can intersect at most $\op{mult}(\ZZ^{n})$ different elements of
the partition $\ZZ^{n}$ for $m_\alpha\le n\le n_\alpha$. We can
assume that $r_\alpha < \epsilon$.

Let $E_n$ denote a $(n,r_\alpha)$-spanning set of minimal
cardinality for $S$ in $X$. For each $x \in X$ consider the
singularity partition $\{Z \cap (\{x\} \times Y)\,|\,Z \in
\ZZ^{n_\alpha(x)}\}$. Each element of this partition is an open
polytope in the fiber $\{x\} \times Y$, and hence we can subdivide
the partition so that each element has diameter no greater than
$\epsilon$. Denote the resulting partition of $\{x\} \times Y$ by
$F_x$. The refinement can be done in such a way that $| F_x| \leq
C_0/\epsilon^d \cdot |P^{n_\alpha}(x)|$ for some constant $C_0 \in
\R_+$.

For $x \in X$ we let $t_0(x)=0$ and define recursively
 $$
t_{k+1}(x)=t_k(x)+|F_{S^{t_k(x)}(x)}|\,.
 $$
Let $q(x)=\min\{k >0 \,|\, t_{k+1}(x) \geq n\}$. We will denote
$q=q(x)$. For $x \in E_n$, $z_0 \in F_x$, $z_1 \in
F_{S^{t_1(x)}(x)},\dots,    z_q \in F_{S^{t_q(x)}(x)}$ denote
\begin{eqnarray*}
V(x;z_0,\dots,z_q)=\big{\{} w \in X\times Y\,\big{|}\, d((S
\tilde\times T)^{t+t_k(x)}(w),(S \tilde\times
T)^{t}(z_k))<2\epsilon
\\
\forall\, 0 \leq t \leq | F_{S^{t_{k-1}(x)}(x)}|,\,0 \leq k \leq
q(x) \big{\}}\,.
\end{eqnarray*}
Then $\cup_{x,z_0,\dots,z_q} V(x;z_0,\dots z_q)=U_n \times Y$ and
for any $(n, 4\epsilon)$-separating set $K \subset X \times Y$ for
$S \tilde\times T$ we have $| K \cap V(x;z_0,\dots z_q)| \leq 1$.
Thus if $K$ is a maximal $(n, 4 \epsilon)$-set, then the
cardinality of $K$ is bounded by the number of ways we can choose
$x,z_0,\dots,z_q$ modulo the partitions specified above. For fixed
$x \in X$ the number $\Pi_x$ of such admissible combinations
satisfies
$$
\Pi_x\leq \prod_{k=0}^{q(x)} | F_{S^{t_k(x)}(x)}|\,.
$$
Since $q(x) \leq n/m_\alpha$ we have:
\begin{eqnarray*}
\log \Pi_x &\leq& (q(x)+1) \log \frac{C_0}{\epsilon^d} +
\sum_{k=0}^{q(x)}\log | P^{n(S^{t_k(x)}(x))}(S^{t_k(x)}(x))|
\\
&\leq& \frac{n+m_\alpha}{m_\alpha} \log \frac{C_0}{\epsilon^d} +
(a+\alpha)\sum_{k=0}^{q(x)} n(S^{t_k(x)}(x))) \\
&\leq& \frac{n+m_\alpha}{m_\alpha} \log \frac{C_0}{\epsilon^d} +
(a+\alpha)(n+n_\alpha) \\
\end{eqnarray*}
Let $Q(S \tilde\times T,n,4 \epsilon)$ denote the cardinality of a
maximal $(n, 4\epsilon)$-separating set for $S \tilde\times T$ in
$X \times Y$. We have that
 $$
\frac1n \log Q(S \tilde\times T,n,4 \epsilon) \leq \frac1n \log
|E_n| + (\frac{1}{m_\alpha}+
\frac1n)\log\frac{C_0}{\epsilon^d}+(a+\alpha)(1+\frac{n_\alpha}{n})\,.
 $$
This yields
$$
\ls \frac{1}{n} \log Q(S \tilde\times T,n,4 \epsilon) \leq \h(S)+
\frac{1}{m_\alpha} \log \frac{C_0}{\epsilon^d}+a+\alpha\,.
$$
Since the left hand side does not depend on $\alpha$ we can let
$\alpha \rightarrow 0$. Then $m_\alpha \rightarrow \infty$ and we
get
$$
\ls \frac{1}{n} \log Q(S \tilde\times T,n,4 \epsilon) \leq
\h(S)+a\,.
$$
Finally let $\epsilon \rightarrow 0$.
 \end{proof}

 \begin{lem} \po \label{lem5}
$\forall {\bf x} \in \Gamma:\,\, \ls \frac1n \log | \Y^{n,{\bf
x}}| \leq \lambda^+(T_{\bf x})+\hm(T_{\bf x})\,.$
 \end{lem}

\noindent The proof of Lemma \ref{lem5} is similar to the proof of
Theorem 1 and will be omitted.

Combining Lemmata \ref{lem4} and \ref{lem5} we obtain:
 \begin{eqnarray*}
\h(S\tilde\times T) &\leq& \h(S)+\sup_{{\bf x}\in \Gamma} \ls
\frac1n \log \Bigl[\op{mult}(S^n,x_0)|\Y^{n,{\bf x}}|\Bigr]  \\
&\leq& \h(S)+\hm(S)+\sup_{{\bf x} \in \Gamma}\ls \frac1n \log  |\Y^{n,{\bf x}}| \\
&\leq& \h(S)+\hm(S)+\sup_{{\bf x} \in \Gamma} \Big{(}
\lambda^+(T_{\bf x})+\hm(T_{\bf x}) \Big{)}\,.
 \end{eqnarray*}

%%%%%%%%%%%%%%%%%%%%%%%%%%%%%%%
%% Proof of Theorem 4
%%%%%%%%%%%%%%%%%%%%%%%%%%%%%%%

\subsection{Proof of Theorem 4}
Let $\mu$ be an $f$-invariant Borel probability measure on $X
\times Y$. Denote the projection of $\mu$ to $X$ by
$\mu_X=\pi_*\mu$ and let $\{\nu_x\}$ be the canonical family of
conditional measures on the fibers $\pi^{-1}(x)$. By the
generalized Abramov-Rokhlin formula \cite{BC} (Bogensch¬tz and
Crauel removed restrictions on the the maps $S$ and $T_x$ in the
original formula \cite{AR}) we have:
$$
h_\mu(S \tilde\times T)=h_{\mu_X}(S)+h_\mu(T|S)\,,
$$
where
$$
h_\mu(T|S,\xi)=\lim_{n \rightarrow \infty} \frac1n \int H_{\nu_x}
\Big{(} \bigvee_{k=0}^{n-1} (T_{S^{k-1}(x)}\circ \dots \circ T_x
)^{-1}(\xi)\Big{)}\,d\mu_X(x)\,,
$$
for a measurable partition $\xi$ of $Y$ and $h_\mu(T|S)=\sup_{\xi}
h_\mu(T|S,\xi)$, where the supremum is taken over all finite
measurable partitions $\xi$ with finite entropy and the refinement
$\bigvee_{k=0}^{n-1} (T_{S^{k-1}(x)}\circ \dots \circ T_x
)^{-1}(\xi)$ is understood with respect to all orbits
$(x,Sx,\dots,S^{k-1}x)$ starting at $x$ (as in \S\ref{S23}).

For $\e>0$ we choose $\xi$ such that $h_\mu(T|S) \leq
h_\mu(T|S,\xi)+\e$. Clearly
\begin{multline*}
\frac1n H_{\nu_x} \Big{(} \bigvee_{k=0}^{n-1} (T_{S^{k-1}(x)}\circ
\dots
\circ T_x )^{-1}(\xi)\Big{)} \leq  \frac1n \log \op{mult}(S^n,x) \\
+\sup_{{\bf x}\in \Gamma_x} \frac1n H_{\nu_x} \Big{(}
\bigvee_{k=0}^{n-1} (T_{(x_{n-1},z_{n-1})}\circ \dots \circ
T_{(x_0,z_0)} )^{-1}(\xi)\Big{)}\,.
\end{multline*}
where $\Gamma_x \subset \Gamma$ is the set of sequences
$((x_0,z_0),(x_1,z_1)\dots)$ with $x_0=x$. If the maps $T_x$ are
non-expanding it follows from the Ruelle-Margulis inequality
\cite{KH} that the last term tends to $h_{\nu_x}(T_{\bf x})=0$ as
$n \rightarrow \infty$. Moreover this happens $\mu_X$-uniformly
and so $h_\mu(T|S)\le\hm(S)+\e$. Let $\e\to0$.

%%%%%%%%%%%%%%%%%%%%%%%%%%%%%%%
%% PROOF OF COROLLARIES
%%%%%%%%%%%%%%%%%%%%%%%%%%%%%%%

\section{Proof of the corollaries}

%%%%%%%%%%%%%%%%%%%%%%%%%%%%%%%
%% PROOF OF COROLLARIES
%%%%%%%%%%%%%%%%%%%%%%%%%%%%%%%

Corollaries \ref{cor1}, \ref{cor2} and \ref{cor7} are obvious.
Corollaries \ref{cor3}, \ref{cor4} follow from the following
section. Corollary \ref{cor5} is implied by the fact that
$\hm(S)=0$ if $\op{dim}(X)=1$ and Corollary \ref{cor8} follows
from our Theorem 4 and Theorem 3 of Buzzi \cite{B1}.

%%%%%%%%%%%%%%%%%%%%%%%%%%%%%%%
%% Estimate for the angle expansion rate
%%%%%%%%%%%%%%%%%%%%%%%%%%%%%%%

\subsection{Estimate for the angles expansion rates}\label{S31.}

Define
 $$
\lambda^+_{[i]}(f)=\mathop{\overline\lim}\limits_{n\to\infty}
\sup_{x\in U_n}\max_{0\le k\le i}\frac1n\log\|\Lambda^kd_xf^n\|.
 $$
Obviously $\lambda^+_{[i]}(f)\le
i\lambda^+_{[1]}(f)=i\cdot\lambda_{\max}(f)$.

 \begin{theorem}\po
The following estimate holds: $\rho_i(f)\le
\lambda^+_{[i]}(f)-i\lambda_{\min}(f)$.
 \end{theorem}

 \begin{proof}
Let us first calculate the differential of the spherical
transformation $^{(s)}\!\!A:S^{d-1}\to S^{d-1}$,
$^{(s)}\!\!A(x)=\dfrac{Ax}{\|Ax\|}$, corresponding to
$A:\R^n\to\R^n$.
 \begin{lem}\po
If $A\in\op{GL}(\R^d)$, then
$d\bigl({}^{(s)}\!\!Ax\bigr)(v)=P^{\perp}_{w}
\biggl(\dfrac{Av}{\|Ax\|}\biggr)$, where $P^{\perp}_{w}$ is the
orthogonal projection along $w={}^{(s)}\!\!A(x)$.
 \end{lem}
In fact, $d\|Ax\|=\langle{}^{(s)}\!\!Ax,d(Ax)\rangle$ and so
 $$
d\bigl({}^{(s)}\!\!Ax\bigr)(v)=\dfrac{Av}{\|Ax\|}-{}^{(s)}\!\!A(x)\cdot
\Big\langle{}^{(s)}\!\!A(x),\dfrac{Av}{\|Ax\|}\Big\rangle.
 $$

From this lemma we get:
$\|d\bigl({}^{(s)}\!\!Ax\bigr)(v)\|\le\dfrac{\|Av\|}{\|Ax\|}\le
C\cdot\dfrac{\max|\op{Sp}(A)|}{\min|\op{Sp}(A)|}$ for
$\|x\|=\|v\|=1$, where the constant $C$ depends on the eigenbasis
of $A$ (in non semi-simple case -- normal basis) only. Since this
eigenbasis is the same for all iterates of $A$ and we have a
finite number of pieces in $f$, we get for all $v$:
$\|d_vd_x^{(s)}\!\!f^n\|\le C_n\cdot
\dfrac{\max|\op{Sp}(d_xf^n)|}{\min|\op{Sp}(d_xf^n)|}$ (with
sub-exponentially growing $C_n$) and so the maximal vertical
(spherical) Lyapunov exponent for the map $d^{(s)}\!f$ at the
point $(x,v)\in STX$ does not exceed the difference
$\overline{\chi}_{\op{max}}(x)-\underline{\chi}_{\op{min}}(x)$.

Similarly, we have:
$\|d\bigl({}^{(s)}\!\!Ax\bigr)\|\le\dfrac{\|\Lambda^kA\|}{\|Ax\|^k}$
and
 $$
\|\Lambda^k d_vd_x^{(s)}\!\!f^n\|\le C_n\cdot
\dfrac{\max\{|\l_1\cdots\l_k|:\l_j\in\op{Sp}(d_xf^n),\l_i\ne\l_j\}}
{\min|\op{Sp}(d_xf^n)|^k}
 $$
with $\ls\frac1n\log C_n=0$ ($\l_i\ne\l_j$ means that the
eigenvalues are different, though in the multiple case they can be
equal), whence the claim.
 \end{proof}
\vspace{4pt}

Note that $\rho_i(f)$ is conformally invariant in the cocycle
sense: The cocycle $\mathcal{A}:X\to\mathop{GL}(\mathbb{R}^n)$,
$x\mapsto d_xf$, can be changed by any cocycle
$\alpha:X\to\mathbb{R}$,
$\mathcal{A}\mapsto\alpha\cdot\mathcal{A}$. Then the Lyapunov-type
characteristics $\rho_i(f)$ do not change.

However if the cocycle $\alpha$ has different upper and lower
Lyapunov exponents, then the quantity
$\lambda_{\max}(f)-\lambda_{\min}(f)$ used in the bound is not
invariant.

%%%%%%%%%%%%%%%%%%%%%%%%%%%%%%%
%%%%%%%%%%%%%%%%%%%%%%%%%%%%%%%

\subsection{Proof of Corollary \ref{cor6}}
It was shown by Buzzi \cite{B1} that $\hm(A)=0$. This was proven
for strictly expanding maps, but the proof extends literally for
any non-degenerate $A$. Thus the bound from above follows from
Theorem 1.

Let's prove that $\h(A)\ge\op{Jac}^+ A$. Assume $A$ is
semi-simple. Let $\l$ be an eigenvalue with $|\l|>1$ and $v$ the
corresponding unit eigenvector. We divide $[0,1]^d$ into domains
$k\e/{|\lambda|^n}\leq \langle x,v \rangle \leq
(k+1)\e/{|\lambda|^n}$. A $(d_n^f, \epsilon)$-ball intersects no
more than two such domains, the total number of which is less than
$\sqrt{d}\cdot|\l|^n/2\e$.

The same holds for other $\l\in\op{Sp}(A)$, so the number of
$(d_n^f,\e)$-balls to cover $[0,1]^d$ is at least
$C_0(\sqrt{d}/2\e)^m(\op{Jac}^+A)^n$, where $m$ is the number of
eigenvalues with absolute value greater than 1 and $C_0$ some
$n$-independent constant.

If $A\in\op{GL}_d(\mathbb{R})$ is not semi-simple, the estimates
change sub-exponentially, implying the same result. Note that the
formula of the theorem holds true even in the case, when $A$ is
degenerate, though arguments should be modified.

%%%%%%%%%%%%%%%%%%%%%%%%%%%%%%%
%% Examples
%%%%%%%%%%%%%%%%%%%%%%%%%%%%%%%

\section{Examples}
\noindent
\begin{figure}
\begin{center}
\includegraphics[width=5.5cm]{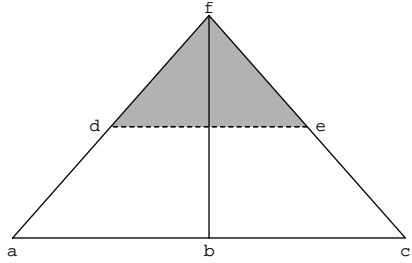}
\caption{{\small Shows the domains of continuity and their image
in Example 1.}} \label{fig1}
\end{center}
\end{figure}
\noindent {\bf Example 1:} Let $X$ be a triangle with vertices in
$a=(-1,0)$, $c=(1,0)$ and $f=(0,1)$. Divide this triangle in two
by taking $X_1$ to be the left triangle with vertices $a=(-1,0)$,
$b=(0,0)$ and $f=(0,1)$. Let $X_2$ be the right triangle with
vertices $b=(0,0)$, $c=(1,0)$ and $f=(0,1)$. Let $X$ be compact,
i.e. the sides are contained in $X$, and let $X_1$ and $X_2$ be
open. Then $\ZZ=\{X_1,X_2\}$ is a finite collection of open
disjoined polytopes in $X$, and $X_1 \cup X_2$ is dense in $X$.

We define a map $S$ on $X'=X_1 \cup X_2$ by the formula
 $$
S(x)=
\begin{cases}
A_1x+B_1 \text{ if } x \in X_1 \\
A_2x+B_2 \text{ if } x \in X_2
\end{cases}\text{where}\,\,A_1=\frac12
\begin{bmatrix} 2 & -1 \\ 0 & 1
\end{bmatrix}, \,\,A_2=\frac12
\begin{bmatrix} 2 & 1 \\ 0 & 1
\end{bmatrix}
 $$
and $B_1=(1/2,1/2)$, $B_2=(-1/2,1/2)$. This maps both $X_1$ and
$X_2$ to the triangle with vertices $d$, $e$ and $f$. Observe that
$S^n(x)$ tends to the point $(0,1)$ with exponential speed for all
$x \in U_f$. This implies that $\h(S)=0$. However the multiplicity
of $\ZZ^n$ at the point $f=(0,1)$ is $2^n$, whence $\hm(S)=\log
2$.

It is also easy to see that $\hs(S) = \log 2$. The eigenvalues of
$A$ are $\frac12, 1$, so the map is non-expanding and
$\lambda^+(f)=0$. Changing $S(x)\mapsto\frac12(S(x)+f)$ we obtain
a strictly contracting piecewise affine map with positive $\hs$
and $\hm$.

The growth of multiplicity is produced by angular expansion: on
$ST_fX\simeq S^1$ the spherization is conjugated to the
map $\theta \mapsto 2\theta$, whence $\rho(S)=\log 2$. \\

\noindent
\begin{figure}
\begin{center}
\includegraphics[width=10.5cm]{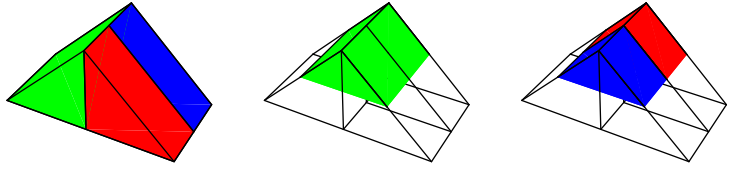}
\caption{{\small The left figure shows the three domains of
continuity for the skew-product $S \tilde\times T$ in Example 2.
Two right figures show the images of the continuity domains.}}
\label{fig2}
\end{center}
\end{figure}

\noindent {\bf Example 2:} The map in Example 1 can be modified to
obtain positive topological entropy. We give here an example of a
piecewise affine non-expanding map with positive topological
entropy in dimension $3$, but it is also possible to construct
such an example in dimension $2$, see \cite{KR2}.

Let $Y=[0,1]$ and let $S:X_1 \cup X_2 \rightarrow X$ be as in
Example 1. For each $x \in X_1 \cup X_2$ we take a piecewise
affine map $T_x \in \op{PAff}(Y,Y)$. For $x \in X_1$ we let
$T_x=Id_Y$ and for $x \in X_2$ we let $T_x$ be the interval
exchange
 $$
T_x(y)=
\begin{cases}
y+1/2 \,\text{ if } y \in (0,1/2), \\
y-1/2 \,\text{ if } y \in (1/2,1).
\end{cases}
 $$
The map $f=S \tilde\times T \in \op{PAff}(X \times Y,X \times Y)$
has three domains of continuity $Z_1,Z_2$ and $Z_3$. These domains
and the images are shown in Figure 2.

As in Example 1 the cardinality of the continuity partition of $S$
grows like $2^n$, but in this example we see that if $x_1,x_2 \in
X$ are elements of different continuity domains for $S$, then
$d(f(x_1,y),f(x_2,y))\geq 1/2$ for all $y\in Y$. Hence the number
of $(d_n^f,\e)$-balls needed to cover $X \times Y$ is at least
$\ZZ^n$, where $\ZZ^n$ is the continuity partition of $f^n$. This
implies $\h(f) \geq \log 2$. The opposite inequality follows from
Theorem 2, whence $\h(f)=\log 2$.

This example has several interesting aspects. First of all it is
an example of a non-expanding piecewise affine map with positive
entropy. We can easily modify $f$ to make it strictly contracting
without changing its topological entropy.

A second important point is that $\h(S)=\h(T_x)=0$ for all $x \in
X'$, so the the entropy of the skew product $S \tilde\times T$
exceeds the combined entropy of its factors. This shows that the
term $\hm(S)$ on the right hand side of the inequality in Theorem
\ref{th3} cannot be removed. Also this justifies Remark \ref{rk5}. \\

\noindent {\bf Example 3:}
\begin{figure}
 \begin{center}
  \begin{picture}(100,100)
 \put(0,0){\circle*{4}}
 \put(100,0){\circle*{4}}
 \put(25,50){\circle*{4}}
 \put(50,50){\circle*{4}}
 \put(75,50){\circle*{4}}
 \put(50,100){\circle*{4}}
 \put(0,0){\line(1,2){50}}
 \put(100,0){\line(-1,2){50}}
 \put(0,0){\line(1,0){100}}
 \put(25,50){\line(1,0){50}}
 \put(50,50){\line(0,1){50}}
 \put(55,98){$A$}
 \put(14,46){$B$}
 \put(47,39){$C$}
 \put(80,46){$D$}
 \put(-10,3){$E$}
 \put(103,3){$F$}
  \end{picture}
 \end{center}
\caption{{\small Example of system with zero Lyapunov exponents
and positive entropy.}} \label{fign}
\end{figure}
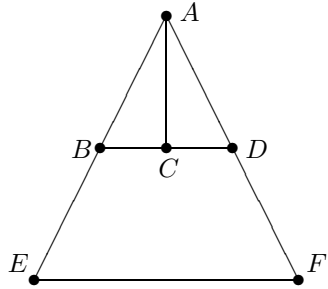
Let us consider a system with continuity domains as in Figure
\ref{fign}. The dynamics $f$ is piece-wise affine and is given by
the rules:
 $$
ABC\longrightarrow AEF,\ \ ADC\longrightarrow AFE,\ \
BDFE\longrightarrow BDFE,
 $$
the last map being the identity. Since every point eventually
comes into the domain $BDFE$, all the Lyapunov exponents vanish.
But the multiplicity of the point $A$ is $\log 2$ and this easily
yields $\h(f)=\log 2$.

This shows that in the estimates of the main theorems we cannot
change the difference $\lambda_{\max}(f)-\lambda_{\min}(f)$ to the
maximal difference of upper and lower Lyapunov exponents
$\sup_{x\in U_f}\max_{i,j}
(\overline{\chi}_i(x)-\underline{\chi}_j(x))$. However we suggest
that the estimate in Theorem \ref{th1} can be refined by changing
$\lambda^+(f)$ to the maximal sum of positive upper Lyapunov
exponents $\sup_{x\in U_f}
\sum\overline{\chi}_i^{\,+}(x)\le\lambda^+(f)$. \\

\noindent {\bf Example 4:}
\begin{figure}
 \begin{center}
  \begin{picture}(100,100)
 \put(0,0){\circle*{4}}
 \put(100,0){\circle*{4}}
 \put(0,100){\circle*{4}}
 \put(100,100){\circle*{4}}
 \put(0,0){\line(0,1){100}}
 \put(0,0){\line(1,0){100}}
 \put(100,100){\line(-1,0){100}}
 \put(100,100){\line(0,-1){100}}
 \put(50,0){\line(0,1){100}}
 \put(1,0){\line(0,1){100}}
 \put(5,0){\line(0,1){100}}
 \put(14,0){\line(0,1){100}}
 \put(18,0){\line(0,1){100}}
 \put(27,0){\line(0,1){100}}
 \put(35,0){\line(0,1){100}}
 \put(40,0){\line(0,1){100}}
 \put(42,0){\line(0,1){100}}
 \put(55,0){\line(0,1){100}}
 \put(58,0){\line(0,1){100}}
 \put(67,0){\line(0,1){100}}
 \put(71,0){\line(0,1){100}}
 \put(76,0){\line(0,1){100}}
 \put(83,0){\line(0,1){100}}
 \put(86,0){\line(0,1){100}}
 \put(89,0){\line(0,1){100}}
 \put(91,0){\line(0,1){100}}
 \put(96,0){\line(0,1){100}}
  \end{picture}
 \end{center}
\caption{{\small A countable piece-wise isometry consisting of
random vertical shifts.}} \label{fign}
\end{figure}
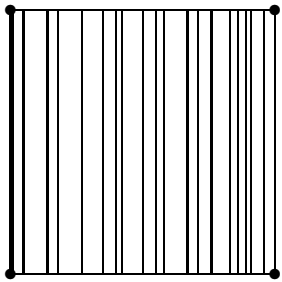
Consider the following map $f$ of $T^2$ to itself. We represent
the torus as a glued square, which is partitioned into countable
number of rectangles $\Pi_i=I_i\times[0,1]$. We define
$f(x,y)=(x+\alpha,y+\beta_i)$ if $x\in I_i$. Thus the map is a
piece-wise isometry with countable number of continuity domains
$\Pi_i$.

The number $\alpha$ is chosen irrational. The intervals $I_i$
(with their lengths $l_i=|I_i|$) and the shift lengths $\beta_i$
are supposed to be sufficiently generic. The value of $\h(f)$
(note that singularity entropy $\hs(f)$ does not have sense here)
depends on the speed of convergence of the series
$\sum_{i=1}^\infty l_i$.

Consider, for instance, the case of rapid convergence, when $l_i$
decrease exponentially or at least polynomially, namely $l_i\le C
i^{-r}$ for some $r>1$ and $C\in\R_+$. Then the frequency with
which an interval of length $\e$ meets a singularity of the base
$\cup\partial I_i\subset S^1$ under rotation by the angle $\alpha$
is at most $p_\e\sim\e^{\frac{r-1}r}$ (for exponential convergence
$p_\e\sim\e\log\frac1\e$). The number of $(d_n^f,\e)$-balls to
cover the torus satisfies: $S(d_n^f,\e)\le
c\cdot(\frac1\e)^{p_\e\cdot(n+2)}$. Consequently $\h(f)=0$.

On the other hand, if the series $\sum_{i=1}^\infty l_i$ converges
slowly, then the entropy may become positive and even infinite.
For example, if $l_i\sim 1/i\log i(\log\log i)^2$, then the
frequency with which an $\e$-interval meets $[1/\e]$ different
intervals $I_i$ under rotation by the angle $\alpha$ has
asymptotic $\sigma_\e\sim 1/\log\log\frac1\e$. Thus choosing the
shifts $\beta_i$ and geometry of the decomposition $I=\cup I_i$
appropriately we may arrive to
$S(d_n^f,\e)\sim(\frac1\e)^{\sigma_\e\cdot(n+2)}$, which yields
$\h(f)=\infty$ (note that $f$ is a skew-product with vanishing
entropies of the base and the fibers).

Note that $\h(f)=0$ for a piecewise isometry $f$ with finite
number of continuity domains \cite{B2} and the same holds for
conformal non-expanding maps \cite{KR2} (and Corollary
\ref{cor2}). For infinite number of
domains this fails. \\

\noindent {\bf Example 5:} In Corollary \ref{cor4} we stated that
if a piecewise affine map $f$ is asymptotically conformal, i.e.
$\lambda_{\op{max}}(f)=\lambda_{\op{min}}(f)$, then $\hm(f)=0$.
This result does not generalize to piecewise smooth maps. This is
shown by the following example due to Buzzi \cite{B2}:

Let $X=[-1,1]^2$ and define $f:X \rightarrow  X$ by
$$
(x,y) \mapsto (x/2,y/2-\op{sgn}(y) x^2)\,.
$$
We observe that
$$
\op{Jac}(f)=\begin{bmatrix} 1/2 & 0 \\ -2x \op{sgn}(y) &
1/2-\delta(y) x^2
\end{bmatrix}\,,
$$
so $\lambda_{\op{max}}(f)=\lambda_{\op{min}}(f)=1/2$. However it
is easy to verify that multiplicity of the origin grows like
$2^n$, whence $\hm(f)=\log 2$.

%%%%%%%%%%%%%%%%%%%%%%%%%%%%%%%
%% refrences
%%%%%%%%%%%%%%%%%%%%%%%%%%%%%%%

\end{document}